\documentclass[12pt]{amsart}
\input{amssymb.sty}
\newtheorem{theorem}{Theorem}[section]
\newtheorem{theorem*}{Theorem}
\newtheorem{lemma}[theorem]{Lemma}

\newtheorem{cor*}{Corollary}

\theoremstyle{definition}

\newtheorem{definition*}{Definition}

\theoremstyle{remark}
\newtheorem{remarks}[theorem]{Remarks}
\newtheorem{remarks*}{Remarks}

\numberwithin{equation}{section}

\newcommand{\Ch}{\operatorname{Ch}}

\newcommand{\Index}{\operatorname{index}}

\newcommand{\rank}{\operatorname{rank}}
\newcommand{\tr}{\operatorname{tr}}

\newcommand{\nc}{\newcommand}
\nc{\hM}{\widehat M}
\nc{\Da}{\Delta}
\nc{\da}{\delta}
\nc{\ta}{\theta}
\nc{\za}{\zeta}
\nc{\A}{\mathcal A}
\nc{\Om}{\Omega}
\nc{\Hj}{{\bar H}_{(2)}^j(\hM)}
\nc{\Hdot}{{\bar H}_{(2)}^{\bullet}(\hM)}
\nc{\Oj}{\Om_{(2)}^j(\hM)}
\nc{\vp}{\varphi}
\nc{\tildvp}{\tilde\vp}
\nc{\hN}{\widehat N}
\nc{\Cj}{C^j(\hK)}
\nc{\HjK}{{\bar H}_{(2)}^j(\hK)}
\nc{\HdotK}{{\bar H}_{(2)}^{\bullet}(\hK)}
\nc{\al}{\alpha}
\nc{\M}{{\mathbb H}}
\nc{\Tau}{T}

\def\al{\alpha}

\def\Om{\Omega}

\def\Todd{{\rm Todd}}
\def\rank{{\rm rank}}

\def\rank{\operatorname{rank}\,}

\def\sup{\operatorname{sup}}
\def\inf{\operatorname{inf}}

\def\C{\mathbb C}
\def\R{\mathbb R}
\def\Z{\mathbb Z}

\def\index{\operatorname{index}}

\def\A{{\mathcal A}}

\def\tr{\operatorname{tr}}

\def\np{\operatorname{\not\!\partial}}

\def\<{\langle}
\def\>{\rangle}

\begin{document}

\hfil{\\}

\title[On positivity of the Kadison constant]
{On positivity of the Kadison constant \\and noncommutative Bloch theory}
\author{Varghese Mathai}
\address{Department of Mathematics, Massachusetts Institute of Technology,
Cambridge, MA 02139, USA}
\email{vmathai@math.mit.edu}
\address{Department of Pure Mathematics, University of Adelaide,
Adelaide 5005, Australia}
\email{vmathai@maths.adelaide.edu.au}
\thanks{This work was completed in part for the Clay Mathematical Institute}
\begin{abstract}
In [V. Mathai, $K$-theory of twisted group
$C^*$-algebras and positive scalar curvature, {\em Contemp.
Math.} {\bf 231} (1999) 203--225],
  we established a natural connection between the
Baum-Connes conjecture and {noncommutative Bloch theory}, viz., the
spectral theory of projectively periodic elliptic operators on covering
spaces. We elaborate on this connection here and
provide significant evidence for a fundamental conjecture in
noncommutative Bloch theory on the non-existence of Cantor set
type spectrum.
This is accomplished by establishing an explicit lower bound
for the Kadison constant
of twisted group $C^*$-algebras in a large number of cases, whenever the
multiplier is rational.
\end{abstract}

\maketitle


\section*{Introduction}
In this paper, we study noncommutative Bloch theory, which is the
spectral theory of elliptic operators on covering
spaces that are invariant under a projective action of the fundamental group.
A central conjecture in noncommutative Bloch theory (cf. \cite{Ma}, \cite{MM})
is the following:\\

\noindent{\bf Conjecture A.} {\em \;
  Suppose that $D : L^2(\widetilde M,
\widetilde E) \to
  L^2(\widetilde M, \widetilde E) $ is a
$(\Gamma,
\bar\sigma)$-invariant,  self-adjoint, elliptic operator on a
normal covering space
$\Gamma\to
\widetilde M \to M$ of a compact manifold $M$,
where $\widetilde E$ is a $\Gamma$-invariant Hermitian vector bundle
over $\widetilde M$. If
$\;\sigma\;$ is rational, then there are only  a finite number of gaps
in the spectrum of $D$ that  lie in any bounded interval $[a, b]$.}\\

Here $\sigma$ is a multiplier on the discrete group $\Gamma$ (i.e., a
normalized $U(1)$-valued group 2-cocycle on $\Gamma$) and
$\bar\sigma$ is its conjugate multiplier. The $(\Gamma,
\bar\sigma)$-action is a projective action of $\Gamma$ that is defined
by the multiplier $\bar\sigma$ (cf. Section 2). A multiplier $\;\sigma\;$
is said to be {\em rational} if
$[\sigma] \in H^2(\Gamma, \mathbb Q/\mathbb Z)$, where we have identified
$\mathbb R/\mathbb Z$ with $U(1)$.

By a fundamental result of \cite{BrSu}, one knows that Conjecture A follows
from the positivity of the Kadison constant $C_\sigma(\Gamma)$
(cf. Subsection 2.3) of the reduced
twisted group $C^*$-algebra $C^*_r(\Gamma, \sigma)$, which is a norm completion
of the twisted group algebra (cf. Section 1). Recall that the Kadison
constant  is the smallest value of the trace evaluated on projections
in $C^*_r(\Gamma, \sigma)$, which can be zero if
there are infinitely many projections
in the algebra. This is the case for instance whenever $\Gamma 
=\mathbb Z^2$ and
$\sigma$ is not a rational multiplier,  \cite{Rieff}.
In this paper, we restrict ourselves to the case when $\Gamma$ is a
torsion-free discrete group and to the case when the
multiplier $\sigma$ has trivial Dixmier-Douady invariant, $\delta(\sigma) =
0$. The method of proving the positivity of the Kadison constant
consists of first defining a twisted analogue of the
Baum-Connes map (or analytic  assembly map or twisted
Kasparov map \cite{Ma}),
\begin{equation}
\mu_\sigma : K_j ( B\Gamma) \rightarrow
K_j (C^*_r(\Gamma, \sigma)),\qquad j=0,1
\end{equation}
which is done in Section 2. It is more convenient for our purposes
than the definition given in \cite{Ma}, but is also equivalent to it.
The {\em twisted Baum-Connes conjecture} asserts that
$\mu_\sigma$ is an isomorphism whenever $\Gamma$ is torsion-free;
the conjecture has to be modified using classifying spaces for proper
$\Gamma$-actions when $\Gamma$ has torsion. Here $B\Gamma$ is the
classifying space of the discrete torsion-free group
$\Gamma$.  In \cite{Ma} it was shown that whenever the discrete group has the
{\em Dirac-Dual Dirac property}, then the twisted Baum-Connes map $\mu_\sigma$
is split injective. For example, by results of Kasparov
\cite{Kas}, discrete  subgroups of connected Lie groups
have the Dirac-Dual Dirac property. If $\Gamma$ is
assumed to be $K$-amenable  in addition, then it has been proved in \cite{Ma}
that $\mu_\sigma$ is an isomorphism. By results of
Kasparov \cite{Kas2},  Julg-Kasparov \cite{JuKas} and Higson-Kasparov
\cite{HiKa}, all discrete amenable groups,  all discrete subgroups of
${\rm \bf SO_0}(1,n)$ and of
${\rm \bf SU_0}(1,n)$ are $K$-amenable and  have the Dirac-Dual Dirac property.

The canonical trace $\;\tau\;$ on the reduced twisted group $C^*$-algebra
$C^*_r(\Gamma, \sigma)$ induces a homomorphism $[\tau] : K_0 (C^*_r(\Gamma,
\sigma))
\to ~\R.\;$ Under the assumption that the twisted
Baum-Connes map
$(1)$ is an isomorphism, we will use a twisted $L^2$
index theorem for covering spaces to compute in Section 2
the range of the trace on $K$-theory in terms of
classical characteristic classes.  This enables us to prove our main
theorem supporting the conjecture  above. We first recall that a $CW$
complex $X$ is said to have finite rational cohomological dimension if
$H^j(X, \mathbb Q) = 0$ for all $j$ large. The cohomological
dimension of $X$ is  then the largest $j$ such that $H^j(X, \mathbb Q)
\ne 0$. A discrete group $\Gamma$ is said to have finite rational
cohomological dimension if its classifying space $B\Gamma$ has
finite rational cohomological dimension.\\

\noindent{\bf Theorem B (Positivity of the Kadison constant).}
{\em \; Suppose that $\Gamma$ is a torsion-free discrete group
satisfying the twisted Baum-Connes conjecture, and that $\Gamma$
has finite rational cohomological dimension. Suppose also that
$\sigma$ is a rational multiplier on $\Gamma$ with trivial 
Dixmier-Douady invariant.
Then the Kadison constant $C_\sigma(\Gamma)$ (cf. Subsection {\rm 2.3}) of
the reduced twisted group $C^*$-algebra $C^*_r(\Gamma, \sigma)$ is
positive.

More precisely, suppose that $\sigma^q = 1$ for some $q> 0$
and $\ell$
is the cohomological dimension of
$\Gamma$, then one has the following lower bound
\begin{equation}
C_\sigma(\Gamma) \ge
c_\ell^{-1} q^{-\ell(\ell + 1)/2}   >0,
\end{equation} where $c_\ell
= \prod_{\ell \ge j > i\ge 0} (2^j - 2^i) >0$.

In particular, Conjecture  {\rm ~A} holds for such pairs $(\Gamma,
\sigma)$.}\\

Observe that by the discussion above,
all discrete subgroups of
${\rm \bf SO_0}(1,n)$ and of
${\rm \bf SU_0}(1,n)$ satisfy the hypotheses of Theorem B,
so do all discrete amenable groups of finite cohomological dimension.
Theorem B was previously known in
the case of two dimensions, when $\Gamma$ is a discrete cocompact subgroup
of ${\rm\bf PSL}(2,\R)$. Here, not only the Kadison constant was 
computed, but also
the range of the trace on $K$-theory was calculated in \cite{CHMM}
when $\Gamma$ is torsion-free (see also \cite{BC}), and when $\Gamma$ has
torsion  in \cite{MM}.
The range of the trace on $K$-theory was computed in the case of 
$\Z^2$ in \cite{Rieff}.
For the relevance of noncommutative Bloch theory
to physics, see \cite{Bel}.

The following theorem is an immediate consequence of Theorem B and the main
theorem in \cite{BrSu} (cf. Theorem \ref{BrSu}),  and is a
quantitative version of the application of
Theorem B to Conjecture A.\\

\noindent{\bf Theorem C.}{\em \; Let $D: L^2(\widetilde M, \widetilde E)
\to L^2(\widetilde M,
\widetilde E)$ be a $(\Gamma, \bar\sigma)$-invariant self-adjoint
elliptic differential operator which is bounded below,
where $\Gamma$ and $\sigma$
satisfy the hypotheses of Theorem {\rm B}. If
$\psi(\lambda)$ denotes the number of gaps in the spectrum
of $D$ that lie in the {\rm 1/2}-line $(-\infty,\lambda]$, then one has
the asymptotic  estimate
\begin{equation}
\lim_{\lambda\to+\infty} \sup\left\{ \frac{(2\pi)^n\;\psi(\lambda)\; 
\lambda^{-n/d}}
{c_\ell q^{-\ell(\ell + 1)/2} \;w_n\; \text{ \rm vol }
(M)} \right\}\le 1,
\end{equation} where $w_n$ denotes the volume of the unit ball in
${\mathbb R}^n$, $\text{ \rm vol }
(M)$ is the volume of $M$,
$n = \dim M$, $d$ is the degree of $D$, $\ell$ is the
cohomological dimension of $B\Gamma$ and $c_\ell
= \prod_{\ell \ge j > i\ge 0} (2^j - 2^i) >0$.}\\

Theorem C gives an asymptotic upper bound for the number of gaps in the
in the spectrum of such operators. In \cite{MS} we prove the existence
of an arbitrarily large number of spectral gaps in the spectrum
of magnetic Schr\"odinger operators with Morse-type potentials on
covering spaces.

Recall that an ICC group is a group such that every nontrivial conjugacy
class in the group is infinite. We also obtain the following theorem in
Section 2.4.\\

\noindent{\bf Theorem D.}{\em \; Suppose that $\Gamma$ and $\sigma$
satisfy the hypotheses of Theorem {\rm B}, and suppose in addition that
$\Gamma$ is an {\rm ICC} group.
Then  there are only a finite number of projections in the twisted group
$C^*$-algebra $C^*_r(\Gamma, \sigma)$, up to Murray-von Neumann
equivalence in the enveloping von Neumann algebra.

In fact, if $\sigma^q = 1$ for some $q>0$ and $\ell$ is the
cohomological dimension of
$B\Gamma$, then there are at most $\; c_\ell q^{{\ell(\ell +
1)}/{2}}\;$ non-trivial projections in the twisted group
$C^*$-algebra $C^*_r(\Gamma, \sigma)$, up to Murray-von Neumann
equivalence in the enveloping von Neumann algebra,
where $c_\ell = \prod_{\ell \ge j > i\ge 0} (2^j - 2^i) >0$.
}\\

We also prove the following useful theorem, which says in particular
that if the Dixmier-Douady invariant $\delta(\sigma) = 0$, then
$K_0 (C^*_r(\Gamma, \sigma)) \cong K_0 (C^*_r(\Gamma)) $, i.e.,
the $K$-theory in this case can be computed from the (untwisted)
Baum-Connes conjecture.\\

\noindent{\bf Theorem E.}{\em \;
Let $\Gamma$ be a discrete group and $\sigma^0, \sigma^1$ be multipliers
on $\Gamma$ such that $\delta(\sigma^0) = \delta(\sigma^1)$. Then there
is an isomorphism
$$
\lambda_1 : K_0(C^*_r(\Gamma, \sigma^1)) \stackrel{\cong}\to
K_0(C^*_r(\Gamma, \sigma^0)).
$$}

Section 1 contains the preliminary material. The twisted Baum-Connes map is
defined in Section 2, and Theorems B, C and D are also
proved there. Theorem E is proved in Section 3.

I would like to thank Stephan Stolz for some very useful remarks on the paper.

\section{Preliminaries}

\subsection{Twisted group $C^*$-algebras}
We begin by reviewing the concept of twisted crossed product
$C^*$ algebras.
Let $A$ be a $C^*$-algebra and $\Gamma$ a discrete group.
Let ${\mbox{Aut}}(A)$ denote the group of $*$-automorphisms of $A$.
Let $\alpha: \Gamma\to {\mbox{Aut}}(A)$ be a representation. Then
$A$ is called a $\Gamma$-$C^*$-algebra. Let $\sigma : \Gamma \times \Gamma
\to {\mbox{U}}(1)$ be a {\em multiplier} satisfying
the following conditions:
\begin{itemize}
\item[(1)] $\alpha_1 = 1$, \ $ \sigma(1,s) = \sigma(s,1) = 1 $ for all
$s \in \Gamma$,
\item[(2)] $\alpha_s\alpha_t =  \alpha_{st}$,
\item[(3)] $\sigma(s,t) \sigma(st, r) = \sigma(s,tr) \sigma(t,r) $ for all
$r, s, t \in \Gamma$.
\end{itemize}

Then $(A, \Gamma, \alpha, \sigma)$ is called a twisted dynamical system.
 From such  a twisted dynamical system, one can construct a {\em twisted crossed
product algebra} as follows.
Let $A(\Gamma,\alpha, \sigma)$ denote the finitely supported, $A$-valued
functions on $\Gamma$ and with multiplication given by
$$
\left(\sum_{\gamma\in \Gamma} a_\gamma 
\gamma\right)\left(\sum_{\gamma'\in \Gamma}
b_{\gamma'} \gamma'\right)
  = \sum_{\gamma\in \Gamma} \sum_{\gamma'\in \Gamma} a_\gamma 
\alpha_\gamma(b_{\gamma'})
\sigma(\gamma,\gamma') \gamma\gamma',
$$
where $a_\gamma, b_{\gamma'} \in A$, and with adjoint or involution given by
$$
\left(\sum_{\gamma\in \Gamma} a_\gamma \gamma\right)^* = 
\sum_{\gamma\in \Gamma}
\alpha_\gamma({\overline a_{\gamma^{-1}}}) 
\sigma(\gamma,\gamma^{-1})^{-1} \gamma .
$$
Recall that a {\em covariant representation} $\pi$ of $(A,\Gamma, 
\alpha, \sigma)$ consists of
a  pair of unitary representations $\pi_A$ of $A$ and $\pi_\Gamma$ of 
$\Gamma$ on a Hilbert
space $\mathcal H$ satisfying
$$
\pi_\Gamma(\gamma)\pi_\Gamma(\gamma') = 
\sigma(\gamma,\gamma')\pi_\Gamma(\gamma\gamma')
$$
for all $\gamma, \gamma' \in \Gamma$ and
$$
\pi_\Gamma(\gamma) \pi_A(a)  \pi_\Gamma(\gamma)^{-1} = \pi_A(\alpha_\gamma(a))
$$
for all $\gamma \in \Gamma, a\in A$. It gives rise to an {\em 
involutive representation}
$\tilde\pi$ of $A(\Gamma, \alpha,\sigma)$ as follows:
$$
\tilde\pi\left(\sum_{\gamma\in \Gamma} a_\gamma \gamma\right) =
\sum_{\gamma\in \Gamma} \pi(a_\gamma) \pi_\Gamma(\gamma).
$$
Define the following norm on $A(\Gamma, \alpha, \sigma)$,
$$
||f|| = \mbox{sup}\left\{ ||\tilde\pi(f)|| ; \pi\;\;\;
\mbox{a covariant representation of}\;\;\; (A,\Gamma, \alpha,
\sigma)\right\}.
$$
The the completion of $A(\Gamma, \alpha,\sigma)$ in this norm is 
called the {\em full twisted
crossed product} and is denoted by $A\rtimes_{\alpha,\sigma} \Gamma$. 
Consider the Hilbert
$\Gamma$-module
$$
\ell^2 (\Gamma, A) = \left\{ \sum_{\gamma\in \Gamma} a_\gamma \gamma ;
\sum_{\gamma\in \Gamma} a_\gamma^*a_\gamma  \;\;\;{\mbox{converges
in}}\;A\right\}.
$$
Then $A(\Gamma, \alpha,\sigma)$ acts on $\ell^2(\Gamma, A)$ and its 
norm completion is called
the  {\em reduced twisted crossed product} and is denoted by 
$A\rtimes_{\alpha,\sigma, r}
\Gamma$. We refer to
\cite{Pe} for more information on these constructions.

We remark that when $\alpha$ is trivial and when $A = \C$, then
$A\rtimes_{\alpha,\sigma,r} \Gamma = C^*_r(\Gamma, \sigma)$ is the reduced
{\em twisted
group} $C^*$-{\em algebra} of $\Gamma$.
The canonical {\em trace} $\tau : C^*_r(\Gamma, \sigma) \to \mathbb C$ is
defined as
$\tau \left( \sum_{\gamma\in \Gamma} a_\gamma \gamma\right) = a_1$.
  Similar remarks apply to the full twisted crossed products.

In this paper, we will only discuss the reduced twisted crossed product
$C^*$-algebras, but many of
our results also apply to the full $C^*$-algebra.

\subsection{Dixmier-Douady invariant} Let $\sigma$ be a multiplier on 
a discrete
group $\Gamma$. Consider the short exact sequence of coefficient groups
\[
    1\to\mathbb{Z} \overset{i}{\to} \mathbb{R}
       \overset{e^{2\pi\sqrt{-1}}}{\longrightarrow} {\rm\bf U}(1) \to 1,
\]
which gives rise to a long exact sequence of cohomology groups (the change of
coefficient groups sequence)
$$
    \cdots \to H^2(\Gamma,\mathbb{Z}) \overset{i_*}{\to}
H^2(\Gamma,\mathbb{R})
       \overset{{e^{2\pi\sqrt{-1}}}_*}{\longrightarrow} H^2(\Gamma, 
{\rm\bf U}(1))
       \overset{\delta}{\to} H^3(\Gamma,\mathbb{Z})
  \overset{i_*}{\to} H^3(\Gamma,\mathbb{R}) \to
\cdots.
$$
Then one definition of the {\em Dixmier-Douady invariant} of $\sigma$ is
$\delta(\sigma) \in H^3(\Gamma,\mathbb{Z})$, \cite{DD}. Notice that 
the Dixmier-Douady
invariant of $\sigma$ is always a torsion element.

\subsection{Classifying space for discrete groups}
Let $\Gamma$ be a discrete group, and $E\Gamma \rightarrow B\Gamma$ be a
locally trivial, principal $\Gamma$-bundle such that $E\Gamma$ is contractible
and $B\Gamma$ is paracompact. Then $B\Gamma$ is unique up to homotopy
and is called the classifying space of $\Gamma$. For example, if $\Gamma$ is a
discrete, torsion-free subgroup of a connected Lie group $G$, then $E\Gamma
= G/K$, where $K$ is a maximal compact subgroup of $G$, and $B\Gamma = \Gamma
\backslash G/K$. For groups with torsion, it turns out (cf. \cite{BCH})
that it is more convenient to consider instead $\underline E\Gamma$,
which is defined  to be the classifying space for proper $\Gamma$-actions
(as opposed to free actions). Then  $\underline E\Gamma$ is unique up to
$\Gamma$-homotopy and coincides with
$ E\Gamma$ when $\Gamma$ is torsion-free. For example, if $\Gamma$ is a
discrete subgroup of a connected Lie group $G$, then $\underline E\Gamma
= G/K$.

\subsection{Topological $K$-homology}
  We now give a brief description of the Baum-Douglas version of
$K$-homology. The basic objects are $K$-cycles. A $K${\em-cycle} on a
topological space is a triple $(M, E, \phi)$, where $M$ is a compact
Spin$^{\mathbb C}$ manifold, $E\to M$ is a complex vecor bundle on $M$, and
$\phi : M \to X$ is a continuous map. Two $K$-cycles $(M, E, \phi)$ and
$(M', E', \phi')$ are said to be
{\em isomorphic} if  there is a diffeomorphism $h: M \to M'$ such that
$h^*(E') \cong E$ and $h^*\phi' = \phi$. Let $\Pi(X)$ denote the
collection of all $K$-cycles on $X$.

\noindent $\bullet$ {\em Bordism}: $(M_i, E_i, \phi_i) \in
\Pi(X)$,
$i=0,1,$ are said to be  {\em bordant} if there is a
triple $(W, E, \phi)$, where $W$ is a compact
Spin$^{\mathbb C}$ manifold with boundary $\partial W$,
$E$ is a complex  vector bundle over $W$ and $\phi : W
\to X$ is a continuous map, such that
$(\partial W, E\big|_{\partial W}, \phi\big|_{\partial W})$ is isomorphic
to the disjoint union $(M_0, E_0, \phi_0) \cup (-M_1, E_1, \phi_1)$. Here
$-M_1$ denotes $M_1$ with the reversed Spin$^{\mathbb C}$ structure.

\noindent $\bullet$ {\em Direct sum}: Suppose that  $(M,
E, \phi) \in \Pi(X)$ and that
$E=E_0\oplus E_1$. Then $(M, E, \phi)$ is isomorphic to $(M, E_0, \phi)
\cup (M, E_1, \phi)$.

\noindent $\bullet$ {\em Vector bundle modification}: Let
$(M, E, \phi) \in \Pi(X)$ and
$H$ be an even dimensional Spin$^{\mathbb C}$
vector bundle over M. Let $\widehat M = S(H\oplus 1)$ denote the sphere
bundle of $H\oplus 1$. Then $\widehat M$ is canonically a  Spin$^{\mathbb C}$
manifold. Let ${\mathcal S}$ denote the bundle of spinors on $H$. Since
$H$ is even dimensional, ${\mathcal S}$ is ${\mathbb Z}_2$-graded,
$$
{\mathcal S} = {\mathcal S}^+ \oplus {\mathcal S}^-
$$
into bundles of $1/2$-spinors on $M$. Define $\widehat E = \pi^*(
{\mathcal S}^{+*} \otimes E)$, where $\pi : \widehat M \to M$ is the
projection. Finally,  $\widehat \phi = \pi^*\phi$. Then $(\widehat M,
\widehat E, \widehat\phi) \in \Pi(X)$ is said to be obtained from
$(M, E, \phi)$ and $H$ by {\em vector bundle modification}.

Let $\;\sim\;$ denote the equivalence relation on
$\Pi(X)$ generated by the operations of bordism, direct
sum and vector bundle modification. Notice that
$\;\sim\;$ preserves the parity of the dimension of
the $K$-cycle. Let $K_0(X)$ denote the quotient
$\Pi^{\rm even}(X)/\sim$, where
$\Pi^{\rm even}(X)$ denotes the set of all even dimensional $K$-cycles in
$\Pi(X)$, and let $K_1(X)$ denote the quotient
$\Pi^{\rm odd}(X)/\sim$, where
$\Pi^{\rm odd}(X)$ denotes the set of all odd dimensional $K$-cycles in
$\Pi(X)$.

\subsection{$K$-theory of $C^*$-algebras}
We briefly recall the definition and some useful facts concerning the
$K$-theory of $C^*$-algebras. Let $A$ be a unital $C^*$-algebra. Recall
that two projections $P$ and $Q$ in $A$ are said to be {\em Murray-von Neumann
equivalent} if there is an element $U\in A$ such that $UU^* = P$ and
$U^*U = Q$. Then the set of Murray-von Neumann equivalence classes
of projections in $A\otimes {\mathcal K}$ is an Abelian semigroup,
where $ {\mathcal K}$ denotes the algebra of compact operators on a
seprable Hilbert space. Then
the associated Grothendieck group is an  Abelian group
which is denoted by  $K_0(A)$. We will make use
of the following lemma in this paper.

\begin{lemma}[\cite{Pe}]\label{ped}
Let $P$ and $Q$ be two projections in a $C^*$-algebra $A$ such that
$||P-Q||<1$. Then they are Murray-von Neumann equivalent. In particular,
suppose that $I\ni t \to P(t)$ is a norm continuous path of projections
in $A$, where $I$ is an interval in $\mathbb R$. Then $P(t)$ are mutually
Murray-von Neumann equivalent.
\end{lemma}

\section{On positivity of the Kadison constant}

\subsection{Twisted Baum-Connes map}
In this section, we will define the twisted Baum-Connes map for an
arbitrary torsion-free discrete group $\Gamma$ and for an arbitrary
multiplier $\sigma$ on $\Gamma$ with trivial Dixmier-Douady invariant.

 From the previous section, elements of $K_0 ( B\Gamma)$ are equivalence
classes of even dimensional $K$-cycles $(M, E, f)$, where $M$ is a compact
Spin$^{\mathbb C}$ manifold, $E\to M$ is a complex vecor bundle on $M$, and
$\phi : M \to B\Gamma$ is a continuous map. Let $\np^{\mathbb C}_E
: L^2(M, {\mathcal S}^+\otimes E) \to L^2(M, {\mathcal S}^-\otimes E)$ denote
the Spin$^{\mathbb C}$  Dirac operator with coefficients in $E$.
Let $\Gamma\to\widetilde{M}\overset{p}{\to}M$ be the covering space of
$M$ such that $\widetilde M = \phi^*(E\Gamma)$. Let
$\widetilde{\np^{\mathbb C}_E}$ be the lift of $\np^{\mathbb C}_E$ to
$\widetilde{M}$,
\[
\widetilde{\np^{\mathbb C}_E}^+\,:\,L^2(\widetilde{M},
\widetilde{\mathcal{S}}^+ \otimes \widetilde{E})\to
L^2(\widetilde{M},\widetilde{\mathcal{S}}^-\otimes \widetilde{E}). \]
Note that $\widetilde{\np^{\mathbb C}_E}$ commutes with the $\Gamma$ action
on $\widetilde{M}$.

Since the Dixmier-Douady invariant $\delta(\sigma) = 0$,
there is an $\R$-valued cohomology class $c$ on
$\Gamma$  such that $e^{2\pi i c} = [
\sigma]$.
Let $\omega$ be a closed 2-form on $M$ such that the cohomology
class of $\omega$ is equal to $\phi^*(c)$,
where $c$ is as before. Note that $\widetilde{\omega}
=p^* \omega=d\eta$ is \emph{exact}, since $E\Gamma$ is contractible. Define
$\nabla=d+\,i\eta$.  Then $\nabla$ is a Hermitian connection on the trivial
line bundle  over $\widetilde{M}$, and the curvature of $\nabla,\
(\nabla)^2=i\,
\widetilde{\omega}$. Then $\nabla$ defines
a projective action of $\Gamma$ on $L^2$ spinors as follows:

First of all, observe that since $\widetilde{\omega}$ is $\Gamma$-invariant,
$0=(\gamma^{-1})^*\widetilde{\omega}-\widetilde{\omega}=
d((\gamma^{-1})^*\eta-\eta)$ for all
$\gamma \in\Gamma$. So $(\gamma^{-1})^*\eta-\eta$ is a closed 1-form on
the  manifold $\widetilde{M}$ whose cohomology class is the pullback of a
cohomology class on $E\Gamma$,
therefore
\[
(\gamma^{-1})^*\eta-\eta=d\phi_\gamma\quad{\text{for
all}}\;\gamma\in\Gamma
\] where $\phi_\gamma$ is a smooth function on $\widetilde{M}$ satisfying
in addition,
$$\phi_\gamma(x_0)=0\quad \text{for some} x_0\in\widetilde{M}\quad
{\text{and all}}\;\gamma
\in\Gamma.$$

Then $\bar\sigma(\gamma,\gamma')=\exp(i\phi_\gamma({\gamma'}^{-1}\cdot x_0))$
defines a multiplier on $\Gamma$, i.e., $\bar\sigma:\Gamma\times\Gamma\to
{\rm\bf U}(1)$  satisfies the following identity for all $\gamma, \gamma',
\gamma''\in \Gamma$:
\[ \bar\sigma(\gamma,\gamma')\bar\sigma(\gamma,\gamma'\gamma'')=
\bar\sigma(\gamma\gamma',\gamma'') \bar\sigma(\gamma',\gamma''). \]
This is verified by observing that
  $\phi_\gamma({\gamma'}^{-1}\cdot x)+\phi_{\gamma'}( 
x)-\phi_{\gamma'\gamma}(x)$
is independent of $x\in\widetilde{M}$ for all $\gamma,\gamma'\in\Gamma$.

For $u \in L^2(\widetilde{M},
\widetilde{\mathcal{S}}\otimes \widetilde{E} )$, let $
S_\gamma u  =  e^{i\phi_\gamma}u$, $\;U_\gamma u =
{\gamma^{-1}}^*u$,
and $T_\gamma=U_\gamma S_\gamma$ be the composition, for all
$\gamma\in \Gamma$.  Then $T$ defines a projective $(\Gamma,
\bar\sigma)$-action on $L^2$-spinors, i.e., \[ T_\gamma
T_{\gamma'}=\bar\sigma(\gamma,\gamma')T_{\gamma\gamma'}.
\]

\begin{lemma}
The twisted Spin$^{\mathbb C}$ Dirac operator on $\widetilde{M}$, \[
\widetilde{\np^{\mathbb C}_E}^+\otimes\nabla\,:\,L^2(\widetilde{M},
\widetilde{\mathcal{S}}^+\otimes \widetilde{E})\to
L^2(\widetilde{M},\widetilde{\mathcal{S}}^-\otimes \widetilde{E}) \]
commutes with the projective $(\Gamma, \bar\sigma)$-action.
\end{lemma}

\begin{proof}
Let $D_\eta = \widetilde{\np^{\mathbb C}_E}^+\otimes\nabla$. Then
$D_\eta = D + ic(\eta)$, where $D= \widetilde{\np^{\mathbb C}_E}^+$
and $c(\eta)$ denotes Clifford multiplication by the one-form
$\eta$.
An easy computation establishes that $U_\gamma D_\eta =
D_{{\gamma^{-1}}^*\eta} U_\gamma$ and that
$S_\gamma D_{{\gamma^{-1}}^*\eta} =   D_\eta S_\gamma
\quad \text{for all}\; \gamma\in
\Gamma$. Then $T_\gamma D_\eta =
D_{\eta} T_\gamma$, where
$T_\gamma = U_\gamma  S_\gamma$ denotes the projective
$(\Gamma,\bar\sigma)$-action.
\end{proof}

Let $D^+$ denote the  twisted Spin$^{\mathbb C}$ Dirac operator on
$\widetilde{M}$, $
\widetilde{\np^{\mathbb C}_E}^+\otimes\nabla$ and
$(D^+)^*=D^-$ its adjoint $
\widetilde{\np^{\mathbb C}_E}^-\otimes\nabla$. Then for
$t>0$,
using the standard Gaussian off-diagonal estimates for the heat kernel
one sees that the heat kernels $e^{-tD^-D^+}$ and $e^{-tD^+D^-}$ are
elements in
$C^*_r(\Gamma, \sigma)\otimes {\mathcal K}$, cf. \cite{BrSu}.
Define the idempotent $e_t(D)\in M_2(C^*_r(\Gamma, \sigma)\otimes
{\mathcal K})\cong C^*_r(\Gamma, \sigma)\otimes {\mathcal K}$ as
follows:
$$
    e_t(D) = \begin{pmatrix} e^{-tD^-D^+} & \displaystyle
e^{-tD^-D^+}\frac{(1-e^{-tD^-D^+})}
    {D^-D^+} D^- \\[+11pt]
e^{-tD^+D^-}{D^+} &
       1- e^{-tD^+D^-} \end{pmatrix}.
$$
It is sometimes known as the Wasserman idempotent. Then the {\em
twisted Baum-Connes map} is
\begin{equation}\label{bc}
\begin{array}{lcl}
\mu_\sigma : K_0 ( B\Gamma) &\longrightarrow & K_0 (C^*_r(\Gamma,\sigma))
\\[+7pt]
\mu_\sigma([M, E, \phi]) &  = &  [e_t(D)] - [E_0]
\end{array}
\end{equation}
where $t>0$ and $E_0$ is the idempotent
$$
E_0 = \begin{pmatrix} 0 & 0 \\ 0 &
       1 \end{pmatrix}.
$$\\

\begin{lemma}
The {twisted Baum-Connes map} is well-defined.
\end{lemma}

\begin{proof}
We need to show that
$
\mu_\sigma([M, E, \phi])
$
defined in Equation $(\ref{bc})$ is independent of $t>0$. By Lemma
\ref{ped}, it suffices to show that
$\;{\mathbb R}^+ \ni t \to e_t(D)$ is a norm continuous family of 
projections in
the $C^*$-algebra $M_2(C^*_r(\Gamma, \sigma)\otimes {\mathcal K})$.
By Duhamel's principle, one has the following identity whenever $t>\epsilon>0$,
$$
e^{-(t\pm\epsilon)D^-D^+} - e^{-tD^-D^+}
= \mp \frac{\epsilon}{t} \int_0^t e^{-sD^-D^+} D^-D^+
e^{-(t-s) (D^-D^+ \pm {\epsilon}/{t}
D^-D^+)} ds.
$$
It follows that there is a constant $C$ independent of $\epsilon$ such that
$$
|| e^{-(t\pm\epsilon)D^-D^+} - e^{-tD^-D^+}|| < C \epsilon.
$$
Using Duhamel's principle repeatedly, one sees
that the family of projections
$\;{\mathbb R}^+ \ni t \to e_t(D)$ is continuous in the norm topology.
\end{proof}

\noindent{\em {\bf Twisted Baum-Connes conjecture}.\; Suppose that
$\Gamma$ is a torsion-free discrete group and that
$\sigma$ is a multiplier on $\Gamma$ with 
trivial Dixmier-Douady invariant. Then the twisted
Baum-Connes map
$$
\mu_\sigma : K_j ( B\Gamma) \rightarrow K_j (C^*_r(\Gamma,
\sigma)),\qquad j=0,1,
$$
is an isomorphism.}\\

It turns out that
the twisted Baum-Connes conjecture is a special case of the  Baum-Connes
conjecture with coefficients (cf. \cite{Ma} for details). The following result
is in \cite{Ma}, where the twisted Baum-Connes map is called the
twisted Kasparov map.

\begin{theorem}[\cite{Ma}] Suppose $\Gamma$ is a torsion-free
discrete group that has the
{Dirac-Dual Dirac property}, and 
$\sigma$ is a multiplier on $\Gamma$ with 
trivial Dixmier-Douady invariant.
Then the twisted Baum-Connes map
$$
\mu_\sigma : K_j ( B\Gamma) \rightarrow K_j (C^*_r(\Gamma,
\sigma)),\quad j=0,1,
$$
is split injective. If in addition, $\Gamma$ is $K$-amenable,
then $\mu_\sigma$ is an isomorphism.
\end{theorem}

For example, by results of Kasparov
\cite{Kas}, discrete  subgroups of connected Lie groups
have the Dirac-Dual Dirac property. By results of
Kasparov \cite{Kas2},  Julg-Kasparov \cite{JuKas}
and Higson-Kasparov \cite{HiKa},
all discrete amenable groups, all discrete subgroups of
${\rm \bf SO_0}(1,n)$ and of
${\rm \bf SU_0}(1,n)$ are $K$-amenable and also
have the Dirac-Dual Dirac property.

\subsection{Characteristic classes}
We recall some basic facts about some well-known
characteristic classes that will be used in this paper,
cf. \cite{Hir}.

Let $E\to M$ be a Hermitian vector bundle over the compact
manifold $M$ that has dimension $n=2m$.
  The {\em Chern classes} of $E$, $c_j(E)$, are by definition
{\em integral } cohomology classes.
The {\em Chern character} of $E$, $\Ch(E)$, is a rational
cohomology class
$$
\Ch(E)=\sum_{r=0}^{m} \Ch_r(E),
$$
where $\Ch_r(E)$ denotes the component of $\Ch(E)$ of degree
$2r$. Then $\Ch_0(E) = \rank(E)$, $\Ch_1(E) = c_1(E)$
and in general
$$
\Ch_r(E) = \frac{1}{r!} P_r(E) \in H^{2r}(M, \mathbb Q),
$$
where $P_r(E) \in  H^{2r}(M, \mathbb Z)$ is a polynomial in the Chern
classes of degree
less than or equal to $r$ with {\em integral} coefficients, that is
determined inductively by the Newton formula
$$
P_r(E) - c_1(E) P_{r-1}(E)\ldots + (-1)^{r-1}c_{r-1}(E) P_1(E)
+(-1)^{r} r c_{r}(E) =0
$$
and by $P_0(E) = \rank(E)$. The next two terms are
$P_1(E) = c_1(E)$, $P_2(E) = c_1(E)^2 - 2 c_2(E)$.

The Todd-genus characteristic class of the Hermitian vector bundle $E$
is a rational cohomology class in $H^{2\bullet}(M,\mathbb Q)$,
$$
\Todd(E) = \sum_{r=0}^{m}\Todd_r(E),
$$
where $\Todd_r(E) $ denotes the component of
$\Todd (E)$ of degree $2r$.  Then $\Todd_r(E)  = B_r Q_r(E)$, where
$Q_r(E)$  is a polynomial in the Chern classes of degree
less than or equal to $r$, with
{\em integral} coefficients, and $B_r \ne 0, B_r \in \mathbb Q$ are the
Bernoulli numbers.

\subsection{Range of the canonical trace} Here we will present some
consequences  of the twisted Baum-Connes conjecture above
and the twisted $L^2$ index theorem for covering spaces (cf. Appendix).

The canonical trace $\;\tau\;$ on $\;{C}^*_r(\Gamma, \sigma)\;$
induces a linear map
$$
[\tau] : K_0 ({C}^*_r(\Gamma, \sigma)) \to \mathbb R,
$$
which is called the {\em trace map} in $K$-theory.
Explicitly, first $\;\tau\;$ extends to matrices with entries in
$\;{C}^*_r(\Gamma, \sigma)\;$ as (with Trace denoting matrix trace):
\[
    \tau(f\otimes r) = {\mbox{Trace}}(r) \tau(f).
\]

Then the extension of $\;\tau\;$ to $K_0$ is given by
$\;[\tau]([e]-[f]) = \tau(e) - \tau(f),\;$ where $e$ and $f$ are
idempotent matrices with entries in $\;{C}^*_r(\Gamma, \sigma)$.
The following result is in \cite{Ma} (see also \cite{BC}).

\begin{theorem}[Range of the trace theorem]\label{range} Suppose that 
$(\Gamma, \sigma)$
satisfies the twisted  Baum-Connes conjecture and
$\delta(\sigma) = 0$.
Then the range of the canonical trace on $K$-theory is
$$
[\tau](K_0(C^*_r(\Gamma, \sigma))) =
\left\{ c_0 \int_M \Todd (M)\wedge e^{\omega}
\wedge\Ch(E) ;\; \mbox{\rm for all}\;(M, E, \phi)
\in \Pi^{\rm even}(B\Gamma)\right\}.
$$

Here
$c_0= {1}/{(2\pi)^{n/2}}$ is a universal constant determined by the
relevant index theorem,
$n = \dim M$, $\Todd$ and $\Ch$ denote the Todd-genus and the Chern character
respectively.
\end{theorem}
\begin{remarks}
The set
$$\left\{ c_0 \int_M \Todd (M)\wedge e^{\omega}
\wedge\Ch(E) ;  \;\mbox{\rm for all}\; (M, E, \phi)
\in \Pi^{\rm even}(B\Gamma)\right\}$$ is a
countable discrete subgroup of $\R$. Note that it is {\em not} in general a
subgroup of $\Z$.
\end{remarks}
\begin{remarks}
When $\Gamma$ is the fundamental group of a compact Riemann surface of positive
genus, it follows from \cite{Rieff} in the genus one case,  and \cite{CHMM}
in the general case, that the set
$$\left\{ c_0 \int_M \Todd (M)\wedge e^{\omega}
\wedge\Ch(E) ;  \;\mbox{\rm for all}\; (M, E, \phi)
\in \Pi^{\rm even}(B\Gamma)\right\}$$
reduces to the countable discrete group $\Z + \theta \Z$, where 
$\theta \in [0,1)$
corresponds to the multiplier $\sigma$ under the isomorphism 
$H^2(\Gamma; {\rm\bf U}(1))
\cong \R/\Z$.
\end{remarks}
\begin{proof}[Proof of Theorem {\rm {\ref{range}}}]

By hypothesis, the
twisted Baum-Connes map is an isomorphism. Therefore to compute the 
range of the
trace map on $K_0(C^*_r(\Gamma, \sigma))$, it suffices to compute the range of
the trace map on elements of the form
$$\mu_\sigma([M, E, \phi]), \qquad [M, E, \phi] \in K_0(B\Gamma).$$
Here $(M, E, \phi)$ is an even parity $K$-cycle over $B\Gamma$.
By the twisted analogue of the $L^2$ index theorem of Atiyah \cite{At} and
Singer \cite{Si} for elliptic operators on a covering space that are invariant
under the projective action of the fundamental group defined by $\sigma$,
and which was stated and used by Gromov \cite{Gr1} (cf. Appendix),
one has
\begin{equation*}
    [\tau](\mu_\sigma([M, E, \phi])) =
       c_0 \int_M \Todd (M)\wedge e^{\omega}
\wedge\Ch(E)
\end{equation*}
as desired.
\end{proof}

\subsection{The Kadison constant and the number of projections in
$C^*_r(\Gamma, \sigma)$}

We begin with a key elementary lemma.

\begin{lemma}\label{elementary}
Suppose that $P(t) = \sum_{j=0}^\ell a_j t^j$ is a
polynomial having the property that there is
$t_0 \in \mathbb N$ for which $P(r t_0) \in \mathbb Z$ for all $r\in 
\mathbb Z$.
Then $P(1) \in c_\ell^{-1} t_0^{-\ell(\ell + 1)/2} \mathbb Z$,
where $c_\ell
= \prod_{\ell \ge j > i\ge 0} (2^j - 2^i).$
\end{lemma}

\begin{proof}
We observe that
$$
\begin{array}{rcl}
a_0 + a_1 t_0 + \ldots + a_\ell t_0^\ell & = & P(t_0) \in \mathbb
Z,\\[+7pt] a_0 + a_1 2 t_0 + \ldots + a_\ell 2^\ell  t_0^\ell & = & P(2
t_0) \in \mathbb Z,\\[+7pt]
\ldots \ldots\ldots\ldots\ldots\ldots\ldots\ldots  & = &
\ldots\ldots\ldots,\\[+7pt] a_0 + a_1 2^\ell t_0 + \ldots + a_\ell
2^{\ell^2} t_0^\ell & = & P(2^\ell t_0) \in \mathbb Z.
\end{array}
$$
That is
$$
T A = Z \in {\mathbb Z}^\ell,
$$
where
$$
T = \left(\begin{array}{lclcl} 1 &  t_0 & t_0^2 & \ldots & t_0^\ell\cr
1 &  2t_0 & 2^2t_0^2 & \ldots & 2^\ell t_0^\ell\cr
\ldots &  \ldots & \ldots & \ldots & \ldots\cr
1 &  2^\ell t_0 & 2^{2 \ell} t_0^2 & \ldots & 2^{\ell^2}t_0^\ell
\end{array}\right)
$$
and
$$
A = \left(\begin{array}{c}
a_0\cr a_1 \cr \vdots \cr a_\ell\end{array}\right),
\qquad
Z = \left(\begin{array}{c}
P(t_0)\cr P(2 t_0) \cr \vdots \cr P(2^\ell t_0)\end{array}\right).
$$
Observe that
$\;{\rm det}(T) = t_0^{{\ell(\ell + 1)}/{2}} {\rm det}(T')\;$, where
$T'$ is the Vandermonde matrix
$$
T' = \left(\begin{array}{lclcl} 1 &  1 & 1 & \ldots & 1\cr
1 &  2 & 2^2 & \ldots & 2^\ell \cr
\ldots &  \ldots & \ldots & \ldots & \ldots\cr
1 &  2^\ell  & 2^{2 \ell}  & \ldots & 2^{\ell^2}
\end{array}\right).
$$
Its determinant is given by the formula
$ {\rm det}(T') = \prod_{\ell \ge j > i\ge 0} (2^j -
2^i) > 0$. Therefore
${\rm det}(T) = t_0^{{\ell(\ell + 1)}/{2}} c_\ell > 0$,
where  $c_\ell
= \prod_{\ell \ge j > i\ge 0} (2^j - 2^i)$. Therefore, $T$ is an invertible
matrix, and one has
$$
A = T^{-1} Z = {\rm det}(T)^{-1} S Z
$$
where $S$ is a matrix with entries in $\mathbb Z$. Therefore
$a_j \in  c_\ell^{-1} t_0^{-\ell(\ell + 1)/2} \mathbb Z$ for
$j=0, \ldots ,\ell$. Therfore
$P(1) = \sum_{j=0}^\ell a_j \in c_\ell^{-1} t_0^{-\ell(\ell + 1)/2} \mathbb Z$
as desired.

\end{proof}

We will now recall the definition of the Kadison constant of a twisted group
$C^*$-algebra.
The {\em Kadison constant} of ${C}^*_r(\Gamma, \sigma)$ is defined by:
$$
C_\sigma(\Gamma) = \inf\{ \tau(P) ; \; P \ \ {\mbox{is a
non-zero projection in}} \ \
{C}^*_r(\Gamma, \sigma) \otimes \mathcal K\}.
$$

\begin{proof}[\bf Proof of Theorem B]
The proof uses Theorem \ref{range} and an analysis of the denominators of the
relevent characteristic classes as discussed in Subsection 2.2. By
assumption, there is a rational cohomology class $c$ such that $e^{2\pi i
c} = [\sigma]$. Therefore there is a  positive integer $q \in \mathbb N$
such that $q[c]$  is an {\em integral} cohomology class.
If $(M, E, \phi)$ is  a $K$-cycle on $B\Gamma$, then the  cohomology class
of $q \omega$ is {\em integral}, where $q[\omega] = q\phi^*[c]$. Note also
that
$rq\omega$ is {\em integral} for all $r\in
\mathbb Z$. It follows that there is a line bundle
$L$ and connection
$\nabla^L$ such that  its curvature $(\nabla^L)^2 = iq\omega$, and also that
the induced connection on the line bundle $L^{\otimes r}$ has curvature
$irq\omega$.
In particular, we see that
$$
c_0 \int_M \Todd (M)\wedge e^{rq\omega} \wedge\Ch(E) = \index(\np^{\mathbb
C}_{E\otimes L^{\otimes r}}) \in \mathbb Z
$$
Observe that since $[\omega] = \phi^*(c)$, it follows that
$[\omega^j] = \phi^*(c\cup c \ldots \cup c) =
  0 \; {\rm if} \; j>\frac{1}{2} \;{\rm cohdim} (B\Gamma) = \ell$, 
where ${\rm cohdim} (B\Gamma)$
denotes the cohomological dimension of $B\Gamma$.
Therefore, for all $r\in \mathbb Z$, one has
$$
\mathbb Z \ni c_0 \int_M \Todd (M)\wedge e^{rq\omega}
\wedge\Ch(E)  = c_0\sum_{j=0}^\ell \frac{r^j q^j}{j!}
\int_M \Todd (M)\wedge \omega^j\wedge\Ch(E)
$$
By Lemma \ref{elementary}, it follows that
$$
c_0 \int_M \Todd (M)\wedge e^{\omega} \wedge\Ch(E)  \in
  c_\ell^{-1} q^{-\ell(\ell + 1)/2} \;\mathbb Z,
$$ where $c_\ell
= \prod_{\ell \ge j > i\ge 0} (2^j - 2^i) >0$.
By Theorem \ref{range}, it follows that the  Kadison constant
$C_\sigma(\Gamma) \ge  c_\ell^{-1} q^{-\ell(\ell + 1)/2} >0$.

\end{proof}

\begin{proof}[\bf Proof of Theorem D]
Let $P$ be a projection in ${C}^*_r(\Gamma, \sigma)$. Then
$1-P$ is also a projection in ${C}^*_r(\Gamma, \sigma)$, and one has
$$
1= \tau(1) = \tau(P) + \tau(1-P).
$$
Each term in the above equation is non-negative. By the
Theorem B and by hypothesis,
it follows that $\tau(P)$ belongs to the finite set
$$\{0, \;c_\ell^{-1} q^{-\ell(\ell + 1)/2}, \;2
c_\ell^{-1} q^{-\ell(\ell + 1)/2},
\ldots,\; 1\}.$$ Since
$\Gamma$ is an ICC group, therefore the enveloping von Neumann
algebra of ${C}^*_r(\Gamma,
\sigma)$ is a factor, cf. \cite{CHMM}.
It follows that for each value of the trace, there is a
unique projection up to Murray-von Neumann equivalence in the
enveloping von Neumann algebra.
Therefore there are at most $c_\ell q^{{\ell(\ell + 1)}/{2}}$
non-trivial  projections in
${C}^*_r(\Gamma, \sigma)$, up to Murray-von Neumann equivalence
in the enveloping von Neumann algebra.
\end{proof}

\subsection{Applications to the spectral theory of projectively periodic
elliptic operators}

In this section, we discuss some quantitative results on the spectrum
of projectively periodic
self-adjoint elliptic operators on covering spaces. In particular,
we formulate a generalization of
the Bethe-Sommerfeld conjecture.

Let $D: L^2(\widetilde M, \widetilde E) \to L^2(\widetilde M,
\widetilde E)$
be a self-adjoint elliptic differential operator that commutes with
the $(\Gamma, \bar\sigma)$-action  defined earlier.
We begin with some basic facts about the spectrum of
such an operator. Recall that the {\em discrete spectrum}
of $D$, ${\rm spec_{disc}}(D)$ consists of all the eigenvalues of $D$ that
have  finite multiplicity, and the {\em essential spectrum} of
$D$, ${\rm spec_{ess}}(D)$ consists of the complement ${\rm spec}(D)
\setminus  {\rm spec_{disc}}(D)$. That is, ${\rm spec_{ess}}(D)$ consists
of the  set of accumulation points of the spectrum of $D$, ${\rm
spec}(D)$. It can be shown that the
{\em discrete spectrum} of $D$ is empty (cf. \cite{MM}).
It follows that ${\rm spec_{ess}}(D) = {\rm spec}(D)$, and in particular
that ${\rm spec_{ess}}(D)$ is unbounded.

Since ${\rm spec}(D)$ is a closed subset of $\R$, its complement
$\R\setminus {\rm spec}(D)$ is an open subset of $\R$, and so it is
the countable union of {\em disjoint} open intervals. Each such interval
is called a {\em gap} in the spectrum of $D$. By the previous
discussion, these also correspond to gaps in ${\rm spec_{ess}}(D)$.
Therefore one can ask the following fundamental question:
  {\em How many gaps are there in the spectrum of $D$?}

This question had been studied previously by Br\"uning and Sunada,
\cite{BrSu}, and we will now discuss some of their results.

Note that in general the spectral projections $E_\lambda =
\chi_{(-\infty, \lambda]}(D)$ of $D$,
do not belong to  ${C}^*_r(\Gamma, \sigma) \otimes \mathcal K$.
However it is a result of \cite{BrSu} that
if $\lambda_0 \not\in {\text{spec}}(D)$,
then $E_{\lambda_0} \in {C}^*_r(\Gamma, \sigma) \otimes \mathcal K$.
The proof uses the off-diagonal decay of the heat kernel
of $D^2$. The following quantitative estimate on the number of
gaps in the spectrum of $D$ in terms of the Kadison constant
$C_\sigma(\Gamma)$ is
the main theorem in \cite{BrSu}.

\begin{theorem}[\cite{BrSu}]\label{BrSu}
Let $D: L^2(\widetilde M, \widetilde E) \to L^2(\widetilde M,
\widetilde E)$ be a $(\Gamma, \bar\sigma)$-invariant self-adjoint
elliptic differential operator which is bounded below. If
$\psi(\lambda)$ denotes the number of gaps in the spectrum
of $D$ that lie in the {\rm 1/2}-line $(-\infty,\lambda]$, then one has
the asymptotic  estimate
\[
\lim_{\lambda\to+\infty} \sup\left\{ \frac{(2\pi)^n\;C_\sigma(\Gamma)
\psi(\lambda)\; \lambda^{-n/d}}
{ \;w_n\; {\rm vol }(M)} \right\}\le 1, \]
where $w_n$ is the volume of the unit ball in ${\mathbb R}^n$,
${\rm{vol} }(M) $ is the volume of $M$,
$n = \dim M$ and $d$ is the degree of $D$.
\end{theorem}

This theorem together with Theorem B immediately establishes
Theorem C.
That is, we have shown that
under the conditions of Theorem B,
the spectrum of a projectively periodic elliptic
operator has countably many gaps which can only accumulate at infinity.
In particular, setting $\bar\sigma = 1$, Theorem C gives
evidence that the following
generalization of the
{\em Bethe-Sommerfeld conjecture} is true, \cite{Ma},
\cite{MM} (see also \cite{KaPe}).\\

\noindent{\bf The Generalized Bethe-Sommerfeld conjecture.}
{\em\; Suppose that $\Gamma$ satisfies the Baum-Connes conjecture.
Then the spectrum of any Hamiltonian $H_{V} = \Delta + V$ on 
$L^2(\widetilde M)$ has
only a {\rm finite} number of gaps, where $\Delta$ denotes the Laplacian
acting on $L^2$ functions on $\widetilde M$ and $V$ is a smooth 
$\Gamma$-invariant
function on $\widetilde M$.}\\

We remark that the Bethe-Sommerfeld conjecture has been
proved by Skriganov \cite{Skri} when $\widetilde M$ is the  Euclidean plane.
See \cite{Gru} for related aspects of noncommutative Bloch theory.

\section{Deformation quantization and the proof of Theorem E}

Theorem E is established using a key result of Rosenberg on the rigidity
of $K$-theory under deformation quantization. We will just sketch the
relevant modifications needed to Rosenberg's argument.

\begin{proof}[\bf Proof of Theorem E]
To prove this, we will define an auxilliary algebra, which is
a  formal deformation quantization of $C^*_r(\Gamma, \sigma^0)$.
Since $\delta(\sigma^0) = \delta(\sigma^1)$, there is an
${\mathbb R}$-valued 2-cocycle $c$ on $\Gamma$ such that
$\sigma^1 = \sigma^0 e^{2\pi i c}$.
Define $C^*_r(\Gamma, \sigma^0)[[\hbar]]$ as being the associative
algebra of formal power series over the ring ${\mathbb C}[[\hbar]]$
of power series over $\mathbb  C$, where the multiplication in
  $C^*_r(\Gamma, \sigma^0)[[\hbar]]$ is defined as
$$
a\star b = ab + \sum_{j=1}^\infty \hbar^j \phi_j(a,b),
$$
where
$$
\phi_j(a,b)(\gamma) = \frac{(2\pi i)^j}{j!}\sum_{\gamma_1\gamma_2 =
\gamma} a(\gamma_1) b(\gamma_2) c^j(\gamma_1, \gamma_2)
$$
are $\mathbb C$-bilinear maps $C^*_r(\Gamma, \sigma^0) \times
C^*_r(\Gamma, \sigma^0) \to C^*_r(\Gamma, \sigma^0)$. Then we observe
that
$$
C^*_r(\Gamma, \sigma^0)[[\hbar]]/(\hbar-s) \cong C^*_r(\Gamma, \sigma^s)
$$
as algebras, where $\;\sigma^s = \sigma^0 e^{2\pi i s c}\;$ and $\;(\hbar
-s)\;$ denotes the ideal generated by $\;\hbar-s, \; s\in{\mathbb R}$.
Therefore the canonical projection map
$$
e_s :  C^*_r(\Gamma, \sigma^0)[[\hbar]] \to C^*_r(\Gamma,\sigma^s)
$$
is a homomorphism of algebras for all $ s \in \mathbb R$. By
adapting the proof  in \cite{Ros}, we see that the induced
map in
$K$-theory is an isomorphism,
$$
(e_s)_* : K_0( C^*_r(\Gamma, \sigma^0)[[\hbar]]) \stackrel{ \cong}{\to}
K_0(C^*_r(\Gamma,\sigma^s)) \qquad  \text{\rm for all}\; s \in \mathbb R.
$$
Then $\lambda_s = (e_0)_* (e_s)_*^{-1}: K_0(C^*_r(\Gamma,\sigma^s))
\to K_0(C^*_r(\Gamma, \sigma^0))$ is an isomorphism for all $s$. Setting
$s=1$, we obtain the desired isomorphism.
\end{proof}

\section*{Appendix: A twisted $L^2$ index theorem for covering spaces}
We use  the notation of the previous section.
  Let $P_+$ and $P_-$ denote the orthogonal projections onto the
nullspace of
$D^+=\widetilde{\np^{\mathbb C}_E}^+\otimes\nabla$ and
$D^-=\widetilde{\np^{\mathbb C}_E}^-\otimes\nabla$,
respectively.  Then $(D^+)^* = D^-$ and one has
\[
D^+ P_+ =
0, \quad {\rm and}\quad D^-
P_- = 0
  \]
By elliptic regularity, it follows that the Schwartz kernels of
$P_+$ and $P_-$ are smooth.
Since the operators $D^+$
and $D^-$
commute with the $(\Gamma,\bar\sigma)$-action,
the same is true for the spectral projections $P+$ and $P_-$.
One can define a semi-finite von Neumann trace on
$(\Gamma,\bar\sigma)$-invariant bounded operators $T$ on the Hilbert
space of
$L^2$-sections of $\widetilde S\otimes \widetilde E$ on the universal
cover, similar to how
Atiyah did in the untwisted case \cite{At},
\[
\tau\left(T\right) =
\int_M\,\tr \left(T(x,x)\right)\,dx,
\] if the Schwartz kernel of $T$ is smooth. It is well defined
since
$$
e^{-i\phi_\gamma(x)} T(\gamma x,\gamma y)\,e^{i\phi_\gamma(y)}
= T(x,y)
$$  for all $x, y \in
\widetilde M$ and $\gamma\in\Gamma$, where we have identified the fiber
at $x$ with the fiber at $\gamma x$ via the isomorphism induced by
$\gamma$.  In particular,
$\tr(T(x,x))$ is a  $\Gamma$-invariant  function on
$\widetilde{M}$.
The $L^2$-index of $D^+$ is by definition
$$
\Index_{L^2}
(D^+)
= \tau_s(P),
$$
where $P = \begin{pmatrix} P_+ & 0 \\ 0 &
       P_- \end{pmatrix}$ and $\tau_s$ denotes the graded trace,
       i.e. the composition of $\tau$ and the
grading operator.
Let $k_+(t,x,y)$ and $k_-(t,x,y)$  denote
the heat kernel of the Dirac operators
$D^-D^+$ and $D^+D^-$, respectively. By general results in
\cite{CGT},
\cite{Roe}, one knows that the heat kernels $k_+(t,x,y)$
and $k_-(t,x,y)$ converge
uniformly over compact subsets of
$\widetilde M \times \widetilde M$ to
$P_+(x,y)$ and $P_-(x,y)$, respectively, as $t \to \infty$.  Therefore
if  $D = \begin{pmatrix} 0 & D^- \\ D^+ &
       0 \end{pmatrix}$, then $e^{-tD^2}=  \begin{pmatrix} e^{-tD^-D^+} &
0
\\ 0 &
       e^{-tD^+D^-} \end{pmatrix}$ and one has
\begin{align*}
  \lim_{t\to\infty}
\tau_s(e^{-tD^2})
& =
\lim_{t\to\infty} \int_M \tr(k_+(t,x,x)) dx
- \lim_{t\to\infty} \int_M \tr(k_-(t,x,x)) dx\\ & = \int_M
\tr(P_+(x,x)) dx - \int_M\tr(P_-(x,x)) dx\\ & = \tau_s(P)\\
& = \Index_{L^2}
(D^+).
\end{align*}
Observe that
\begin{align*}
\frac{\partial}{\partial t}
\tau_s( e^{-tD^2})
& = - \tau_s( D^2
e^{-tD^2})\\
& = -  \tau_s( [D,
  D e^{-tD^2}])\\
& = 0,
\end{align*}
since $D$ is an odd operator.
Therefore one also has the analogue of the McKean-Singer formula
in this context, that is, for $t>0$,
$$
\Index_{L^2}
(D^+) = \tau_s(e^{-tD^2})
=  \lim_{t\to\infty}
\tau_s(e^{-tD^2}).
$$
By the discussion below Lemma 2.1, we have
$$
  \tau(\mu_\sigma([M, E, \phi])) = \tau(e_t(D)) - \tau(E_0) =
\tau_s(e^{-tD^2}).
$$
Combining the observations above, one has
$$
  \tau(\mu_\sigma([M, E, \phi])) = \Index_{L^2}
(D^+).
$$
The following twisted analogue of the $L^2$ index theorem for covering spaces
due to Atiyah \cite{At} and Singer \cite{Si}, was stated in 
\cite{Gr1} and proved
for instance in the appendix of \cite{Ma}. It is established using the local
index index theorem (cf. \cite{Get}) and the
observations above.

\begin{theorem}[$L^2$ Index Theorem for $(\Gamma,\bar\sigma )$-invariant Dirac
type operators] The $L^2$ index theorem for elliptic operators which are of
Dirac type  and which are invariant under the projective
$(\Gamma,\bar\sigma)$-action is
$$\Index_{L^2}
(\widetilde{\np^{\mathbb C}_E}^+\otimes\nabla  )
= \int_M \Todd(\Omega)\, e^{\omega}\, \tr(e^{R^E}),
$$
where $\Todd(\Omega)$ denotes the Todd genus of the
Spin$^{\mathbb C}$ manifold $M$ and $\tr(e^{R^E})$  is
the Chern character of the Hermitian vector bundle $E$
over $M$.
\end{theorem}

\end{document}